\magnification =1100
 \baselineskip =13pt
\overfullrule =0pt
\input amstex
\documentstyle{amsppt}
\nologo
 
\topmatter

\title The Gauss map and a noncompact Riemann-Roch formula
for constructible sheaves on semiabelian  varieties
\endtitle

\rightheadtext{Noncompact Riemann-Roch}

 \author{ J. Franecki, M. Kapranov}\endauthor
\address{J.F.: Department of Mathematics, Loyola University, 6525 N. Sheridan Rd.
Chicago IL 60620,\hfill\break 
M.K.: Department of Mathematics, Northwestern  University,
Evanston IL 60208}\endaddress

\email {franecki\@math.luc.edu, kapranov\@math.nwu.edu}\endemail

\endtopmatter

\document

{\bf Abstract:} {For an irreducible subvariety $Z$ in an algebraic group $G$ we define an
integer $\text{gdeg}(Z)\geq 0$ as the degree, in a certain sense, of the
Gauss map of $Z$. It can be regarded as a substitution for the
intersection index of the conormal bundle to $Z$ with the zero
section of $T^*G$, even though $G$ may not be compact. For $G$ a
semiabelian variety (in particular, an algebraic torus $(\Bbb C^*)^n$), we
prove a Riemann-Roch type formula for constructible sheaves on $G$
which involves our substitutions for the intersection indicies. As a
corollary, we get that a perverse sheaf on such a $G$ has 
nonnegative Euler characteristic, generalizing a theorem of
Loeser-Sabbah. 
}

\vskip .3cm

\heading{1. Introduction and statement of results} 
\endheading

Let $X$ be a smooth algebraic variety over $\Bbb C$ and $\frak{F}$ be a
constructible sheaf of $\Bbb C$-vector spaces on $X$. As in other
situations, we have the Riemann-Roch problem: express
$\chi(X,\frak{F}) = \sum_i (-1)^i \dim H^i(X,\frak{F})$ in
terms of some intrinsic geometric invariants of $X$ and
$\frak{F}$. One such invariant is the characteristic cycle
$CC(\frak{F})$ which is a formal $\Bbb Z$-linear combination
${\sum}_\nu n_\nu [\Lambda_{\nu}]$ of irreducible conic
Lagrangian subvarieties $\Lambda_\nu$ in the cotangent bundle
$T^*X$, see [14] for background. When $X$ is compact, the
Riemann-Roch problem has a nice solution, namely [13]:

$$\chi(X,\frak{F}) = ( CC(\frak{F}), [X])_{T^*X}, \leqno (1.1)
$$
where the right hand side is the intersection index, in $T^*X$, of
$CC(\frak{F})$ and the zero section $X \subset T^*X$. It can be
calculated, for instance, by deforming $X$ to the graph of a
$C^{\infty}$ 1-form so that the intersection becomes transverse, and
then counting intersection points (with multiplicities and signs). 

Both the definition of $CC(\frak{F})$ and the formula (1.1) (for
compact $X$) extend to the case when $\frak{F}$ is a bounded
constructible complex (i.e., a complex of sheaves with constructible
cohomology). 

When $X$ is not compact, $\chi(X,\frak{F})$ still makes sense, but
(1.1) is not applicable. We face, therefore, an interesting
{\it noncompact Riemann-Roch problem} of finding $\chi(X, \frak{F})$ in
terms of invariants intrinsic to $X$ (in particular, not involving the
choice of compactification).

The purpose of this paper is to exhibit such a ``noncompact
Riemann-Roch formula'' in a particular class of situations. Namely,
suppose that $X=G$ is an algebraic group with Lie algebra
$\frak{g}$. For $\gamma \in \frak{g}^*$, let $\omega_{\gamma}$ be the
corresponding left-invariant 1-form on $G$, and $\Omega_{\gamma}
\subset T^*G$ be its graph. The $\Omega_{\gamma}$ then form a natural
family of deformations of $X$ and we can use them to make sense of the
intersection index in (1.1) even when $G$ is not compact. More
precisely, if $\Lambda \subset T^*G$ is an irreducible conic
Lagrangian subvariety and $\gamma \in \frak{g}^*$ is {\it generic},
then $\Lambda \cap \Omega_{\gamma}$ consists of finitely many
transversal intersection points; their number will be denoted
$\text{gdeg}(\Lambda)$ and called the {\it Gaussian degree} of $\Lambda$. To
explain the name, recall that $\Lambda$ has the form $T_Z^*X$ for an
irreducible subvariety $Z \subset X$ (notation: $T_Z^*X$ will always
mean the closure of the conormal bundle to the smooth locus $Z_{sm}$ of
$Z$). Denoting $k = \dim(Z)$, we have the left Gauss map
$$\Gamma_Z: Z \to G(k,\frak{g}) \; , \; z \mapsto z^{-1}(T_zZ)\subset T_eG
= \frak{g}$$
which is a rational map, regular on $Z_{sm}$. The number
$\text{gdeg}(\Lambda)$ is the degree of $\Gamma_Z$ in an appropriate sense,
see section 2. For example, if $Z$ is a hypersurface, then the source
and the target of $\Gamma_Z$ have the same dimension, and
$\text{gdeg}(\Lambda)$ is the degree of $\Gamma_Z$ in the usual sense.

Note that $\text{gdeg}(\Lambda) \geq 0$ by construction. We now formulate the
main result of this note. Recall [2] that a semiabelian variety is
an algebraic group $G$ which is an extension
$$1 \to T \to G \to A \to 1, \leqno (1.2)$$
where $A$ is an abelian variety and $T \cong (\Bbb C^*)^n$ is an algebraic
torus.

\proclaim{(1.3) Theorem}
Let $G$ be a semiabelian variety, $\frak{F}$ a bounded constructible
complex on $G$ and $CC(\frak{F})=\sum_\nu
n_{\nu}[\Lambda_{\nu}]$. Then 
$$
\chi(G,\frak{F}) = \sum_{\nu} n_{\nu} \roman{gdeg}(\Lambda_{\nu}).
$$
\endproclaim

\proclaim {(1.4) Corollary}
If, in the situation of (1.3), $\frak{F}$ is a perverse sheaf, then
$\chi(G,\frak{F}) \geq 0$.
\endproclaim

Indeed, for $\frak{F}$ a perverse sheaf, all $n_{\nu} \geq 0$. Here we
use the conventions of [1] for the definition of (middle)
perversity.

\proclaim{(1.5) Corollary}
If $G$ is semiabelian and $Z \subset G$ is a smooth closed subvariety,
then the number $(-1)^{\dim(Z)} \chi(Z,\Bbb C)$ is nonnegative and coincides with
$\roman{gdeg}(Z)$.
\endproclaim

Indeed, $\underline{\Bbb C}_Z[\dim(Z)]$ is perverse and its
characteristic cycle is $[T^*_Z G]$ taken with multiplicity 1.

Corollary 1.4 for $G = (\Bbb C^*)^n$ was proven by Loeser-Sabbah [12]
and given a different proof (applicable to etale sheaves) by
Gabber-Loeser [4]. In the case when $G$ is an abelian variety,
Corollary 1.4 seems to be new, even though Theorem 1.3 in this case
does not need a special proof, being a consequence of (1.1). Moreover,
since the general proof of Theorem 1.3 given below is extremely
simple and transparent, we believe that our approach exhibits the true
reason behind the Loeser-Sabbah observation. Further, another result
of [4] identifies irreducible perverse sheaves $\frak{F}$ on
$(\Bbb C^*)^n$ such that $\chi((\Bbb C^*)^n,\frak{F}) = 1$, with complexes of
solutions of hypergeometric systems (essentially, of the
$A$-hypergeometric systems [5][11], see [6] for a
comparison of the two points of view). On the other hand, in [12]
the second author classified irreducible hypersurfaces in $(\Bbb C^*)^n$
for which the Gauss map has degree 1 and identified them with
(reduced) $A$-discriminantal hypersurfaces. The latter describe the
characteristic varieties of the $A$-hypergeometric system. Thus our
approach explains the analogy between these two results.

Unfortunately, Theorem 1.3 cannot be straightforwardly generalized to
more general algebraic groups. For example, if $G$ is nonabelian
reductive, then it contains affine spaces $\Bbb A^m$ for $m$ ranging from
0 to the number of positive roots of $G$. The Euler characteristic of
$\Bbb A^m$ being 1, Corollary 1.5 (and thus Theorem 1.3)
cannot hold. The same applies when $G$,
while commutative, contains $G_a$, the additive group. Nevertheless,
we believe that there should exist generalizations of Theorem 1.3
which involve some particular classes of constructible sheaves and
complexes.

We are grateful to V. Ginzburg for several useful remarks on the first
version of the text. 
A part of the results of this paper were included in the thesis of the
first author [3]. The second author was partially supported by the
National Science Foundation.

\heading{2. The Gauss map and the Gaussian degree}
\endheading

Let $G$ be a complex algebraic group, $\frak{g}$ its Lie algebra and
$G(k,\frak{g})$ the Grassmannian of $k$-dimensional linear subspaces
in $\frak{g}$. If $Z \subset G$ is an irreducible $k$-dimensional
subvariety, we have the rational map
$$
\Gamma_Z: Z \to G(k,\frak{g})
$$
called the (left) Gauss map and defined as follows. For $x \in G$ let
$$
l_x: G \to G \; , \; l_x(y) = xy
$$
be the left translation by $x$. Then, for a smooth point $z \in Z$ the
value $\Gamma_Z(z)$ is the image of 
$$
d_xl_{x^{-1}}: T_zZ \to T_eG = \frak{g}.
$$
We want to associate to $\Gamma_Z$ a nonnegative integer called its
degree. To this end, consider first a more general situation.

Let $M$ be an irreducible $k$-dimensional variety, $V$ an
$n$-dimensional vector space and $f: M \to G(k,V)$ a rational
map. Replacing, if necessary, $M$ with its Zariski open subset, we can
assume that $f$ is regular. Consider the flag variety $F(k,n-1,V)$ and
its projections
$$G(k,V) \buildrel{p}\over\longleftarrow F(k,n-1,V)
\buildrel {q}\over\longrightarrow G(n-1,V) = P(V^*)$$
Let $\widetilde{M} ={G(k,V)}{\times}_M F(k,n-1,V)$ be the fiber product
with respect to $f$ and $p$. Since $p$ is a smooth map with
$(n-k-1)$-dimensional fibers, $\widetilde{M}$ is an irreducible variety of
dimension $n-1$. The map $q$ induces a regular map 
$q': \widetilde{M} \to \Bbb P(V^*)$ whose source and target have the same dimension. We define
$\deg(f)$ to be the degree of the map $q'$. The following is then
clear.

\proclaim {(2.1) Proposition}
If $W \subset V$ is a generic hyperplane, then $\deg(f)$ is equal to
the number of $x \in M$ such that $f(x) \in G(k,W) \subset G(k,V)$. For
any such $x$ the map $f$ is locally (in the analytic or etale
cohomology) an embedding near $x$ and the intersection of $f(M)$ and
$G(k,W)$ is transversal at $x$.
\endproclaim

When $k=n-1$, we have $\widetilde{M}=M$ and $\deg(f)$ is the degree of $f$
in the usual sense.

We now specialize to the case $M=Z$, $V=\frak{g}$, $n = \dim(G)$ and $f
= \Gamma_Z$. The number $\deg(\Gamma_Z)$ will be denoted by $\text{gdeg}(Z)$
and called the Gaussian degree of $Z$.

Let $\Lambda = T_Z^*G$ be the conic Lagrangian variety associated to
$Z$. We will write $\text{gdeg}(\Lambda)$ for $\text{gdeg}(Z)$. As in
 Section 1, for
$\gamma \in \frak{g}^*$ let $\Omega_{\gamma} \subset T^*X$ be the
graph of the left-invariant1-form $\omega_{\gamma}$ on $G$ associated
to $\gamma$. Proposition 2.1 implies easily:

\proclaim {(2.2) Proposition} Let $\gamma \in \frak{g}^*$ be a generic linear functional. Then
$\Lambda \cap \Omega_{\gamma}$ consists of finitely many points which
are smooth on $\Lambda$ and in which the intersection is
transverse. The number of these points is equal to $\roman{gdeg}(\Lambda)$.
\endproclaim

\heading{3. Characteristic cycle of an open embedding}
\endheading

A non-intrinsic way to find $\chi(U,\frak{F})$ where $U$ is a
noncompact manifold, is to apply (1.1) to $Rj_*\frak{F}$ where $j: U
\hookrightarrow X$ is a smooth compactification. We will indeed use this
approach in the proof of Theorem 1.3, so we recall the (now well
known) proceedure of finding $CC(Rj_*\frak{F})$ from $CC(\frak{F})$,
see [7] [17].

Let $X$ be a not necessarily compact smooth variety, $f \in \Bbb C[X]$ a
regular function, $U \subset X$ the open set $\{ f \neq 0 \}$, and $j:
U \hookrightarrow X$ the embedding. Let $\Lambda \subset T^*U$ be an
irreducible conic Lagrangian variety. For $s \in \Bbb C^*$ let
$$
\Lambda_s^{\#} = \Lambda + s d\log f = \{(\xi + s(d\log f)(x),x) |
(x,\xi) \in \Lambda \}\leqno (3.1)$$
This is a closed (no longer conic) Lagrangian subvariety in
$T^*X$. The total space of the family of $\Lambda_s^{\#}$ is a
subvariety $\Lambda^{\#} \subset T^*X \times \Bbb C^*$. The limit
$\lim\limits_{s\to 0} \Lambda_s^{\#}$ (also called the
specialization of $\Lambda^{\#}$ in [7]) is an effective
Lagrangian cycle in $T^*X$ defined as follows. We first take the
closure $\overline{\Lambda^{\#}}$ in $T^*X \times \Bbb C$ and then form
the scheme-theoretic intersection $\overline{\Lambda^{\#}} \cap (T^*X
\times \{0\})$. The cycle $\lim\limits_{s\to 0} \Lambda_s^{\#}$ is
obtained by taking the irreducible components of this intersection
with the multiplicities given by the scheme structure.

We extend this construction by $\Bbb Z$-linearity to conic Lagrangian
cycles in $T^*U$. Thus, if $\Sigma$ is such a cycle, we have the family
of non-conic cycles $\Sigma_s^{\#}$, $s \in \Bbb C^*$ and the conic cycle
$\lim\limits_{s\to 0} \Sigma_s^{\#}$ in $T^*X$. Now, the fact we
need  is as
follows.

\proclaim {(3.2) Theorem}
If $\frak{F}$ is a bounded constructible complex on $U$, then
$$CC(Rj_*\frak{F}) = \lim_{s\to 0} \, CC(\frak{F})_s^{\#}.$$
\endproclaim

This statement can be obtained from Theorem 3.2 of [7] by applying the
Riemann-Hilbert correspondence, or from Theorem 3.1 of [17]
which is applicable to the more general case of $\Bbb R$-constructible
sheaves.
(To be precise, the concepts of the characteristic cycle used in
[7-9] and [13][14][17] refer to different contexts: holonomic D-modules
vs. constructible complexes. The compatibility of these two definitions
of the characteristic cycle inder the  Riemann-Hilbert correspondence
follows from the results of [7].)

\vskip .2cm 

Consider now a nominally more general situation (cf. [9], Appendix A):
 let $X$ be as before
but suppose that we have $n$ regular functions $f_1, \dots, f_n \in
\Bbb C[X]$. Let $U$ be the intersection of the $n$ open sets $\{f_i \neq
0\}$ and $j: U \hookrightarrow X$ be the embedding. Of course, this
situation can be analyzed by applying Theorem 3.2 to $f = f_1\dots
f_n$, but it will be convenient for us to have a more flexible
formulation.

For a point $s = (s_1, \dots, s_n) \in (\Bbb C^*)^n$ and a conic
Lagrangian variety $\Lambda \subset T^*U$ we form, similarly to (3.1),
a non-conic Lagrangian variety
$$\Lambda_s^{\#} = \Lambda + \sum_{i=1}^n s_i d\log f_i \subset T^*X.
\leqno (3.3)$$
The total space of this family lies in $T^*X \times (\Bbb C^*)^n$. Taking
the closure in $T^*X \times \Bbb C^n$ and then intersecting with $T^*X
\times \{(0, \dots, 0)\}$ defines, similarly to the above, a conic
Lagrangian cycle $\lim\limits_{s \to (0,\dots,0)}
\Lambda_s^{\#}$. Of course, this ``limit'' could be taken along any
curve in $\Bbb C^n$ passing through 0 and generically lying in
$(\Bbb C^*)^n$. As before, we extend this construction by linearity to
conic Lagrangian cycles in $T^*U$. The next statement follows by
iterated application of Theorem 3.2.

\proclaim {(3.4) Theorem}
If $\frak{F}$ is a bounded constructible complex on $U$, then
$$CC(Rj_*\frak{F}) =  \lim_{s \to (0,\dots,0)}
CC(\frak{F})_s^{\#} \; , \; s=(s_1,\dots,s_n) \in (\Bbb C^*)^n.$$
\endproclaim

Taking the limit along different curves approaching $(0,\dots,0)$
corresponds, roughly, to different choices of our equation for the
reducible hypersurface $X \setminus U$. For example, restricting to
the curve with parametric equation $s_i = t^{m_i}$, $m_i > 0$
corresponds to taking $\prod_i f_i^{m_i}$ as an equation.

\vskip .2cm

We now need a slight globalization of Theorem 3.4. First of all, let $(L, \nabla)$
be a line bundle on $X$ with an algebraic flat connection. If $f$ is a
regular section of $L$, then $f^{-1}\nabla f$ is a scalar 1-form regular
over the open set $\{f \neq 0 \}$. We denote this form $\nabla \log
f$.

Suppose now that we have $n$ line bundles with flat connections $(L_i,
\nabla_i)$ on $X$, $i=1,\dots,n$. Suppose $f_i \in \Gamma(X,L_i)$,
$i=1,\dots,n$ and $U \subset X$ is the intersection of the open sets
$\{f_i \neq 0 \}$. As before, let $j: U \hookrightarrow X$ be the
embedding.

\proclaim{(3.5) Theorem}
For a bounded constructible complex $\frak{F}$ on $U$ we have
$$CC(Rj_*\frak{F}) = \lim_{s \to (0,\dots,0)} \biggl(CC(\frak{F}) +
\sum_{i=1}^n s_i \nabla_i \log f_i\biggr).$$
\endproclaim

\noindent {\sl Proof: }
As before, it is enough to consider the case $n=1$, as the
general case can be obtained by iteration, as in [9], Appendix A.
The statement for $n=1$ is a consequence of Theorem 6.3 of [8]
which deals with the more general case of the zero locus of a section
$f$ of an arbitrary line bundle $L$, not necessarily with
connection. The recipe in this case is to consider the ``twisted
cotangent bundles" $(T^*X)^{(s)}, s\in \Bbb C$, defined as the
symplectic quotients of $T^*L$ by the hamiltonian action of $\Bbb C^*$
induced by dilations of $L$. Now, if $L$ is equipped with a flat
connection $\nabla$, then all the $(T^*X)^{(s)}$ become
identified with $T^*X$ and the formulation of [8], Theorem 6.3
reduces to our statement.

\heading {4. Proof of Theorem 1.3 for $G = (\Bbb C^*)^n$}
\endheading

We first consider the case when $G = (\Bbb C^*)^n$ is an algebraic
torus. Let $z_1,\dots,z_n$, $z_i \neq 0$, be the standard coordinates
in $(\Bbb C^*)^n$. We compactify $G$ by the projective space $\Bbb P\,^n$
with homogeneous coordinates $(t_0:\dots:t_n)$ by
$$
j:(\Bbb C^*)^n \hookrightarrow \Bbb P\,^n \; , \; (z_1,\dots,z_n) \mapsto
(1:z_1:\dots:z_n)
$$
For $\nu = 0,\dots,n$ let $A_{\nu}^n \subset \Bbb P\,^n$ be the affine
chart given by $t_{\nu} \neq 0$. This is an affine space with
coordinates $z_i^{(\nu)}$, $i \in \{0,\dots,n\} \setminus \{\nu\}$
given by $z_i^{(\nu)} = \frac{t_i}{t_{\nu}}$. Denote by
$$(\Bbb C^*)^n \buildrel{j_{\nu}}\over\hookrightarrow A_{\nu}^n
\buildrel {k_{\nu}}\over\hookrightarrow \Bbb P\,^n$$
the embeddings. For $\nu=0$ we have $z_i^{(0)} = z_i$.

We now apply Theorem 3.4 to $U=(\Bbb C^*)^n$, $X=A_0^n$, $f_i=z_i$, and
our constructible complex $\frak{F}$. The recipe of the theorem
requires us to introduce the family of 1-forms
$$
\omega_s = \sum_{i=1}^n s_i \, d\log z_i \; , \; s=(s_1,\dots,s_n) \in
(\Bbb C^*)^n.
$$
These forms are precisely the invariant 1-forms on $G=(\Bbb C^*)^n$. We
can view $s$ as an element of $\frak{g}^*$, where $\frak{g} = \Bbb C^n$ is
the Lie algebra of $G$. Theorem 3.4 then gives us:
$$CC(Rj_{0*}\frak{F}) = CC(Rj_*\frak{F})|_{A_0^n} =\lim
_{s\to 0} (CC(\frak{F}) + \omega_s),\leqno (4.1)$$
the limit being taken in $T^*A_0^n$.

Next, we apply Theorem 3.4 to $U = (\Bbb C^*)^n$ and $X = A_{\nu}^n$ with
arbitrary $\nu \in \{0,\dots,n\}$. Then we should take $f_i =
z_i^{(\nu)}$, $i \in \{0,\dots,n\} \setminus \{ \nu \}$ and consider
the 1-forms
$$
\omega_{s'}^{(\nu)} = \sum_{i \neq \nu} s'_i \, d\log z^{(\nu)}_i \; ,
\; s' \in (\Bbb C^*)^{\{0,\dots,n\} \setminus \{ \nu \} }.
$$
Now, each $z^{(\nu)}_i$ is a Laurent monomial in the $z_1,\dots,z_n$,
so $d\log z^{(\nu)}_i$ is an (integer) linear combination of $d\log
z_1, \dots, d\log z_n$. Therefore, $\omega^{(\nu)}_{s'} = \omega_s$
where $s=\phi_{\nu}(s')$ is an image of $s$ under a linear
transformation $\phi_{\nu}: \Bbb C^{\{0,\dots,n\} \setminus \{\nu\} } \to
\Bbb C^n$. It is clear that $\phi_{\nu}$ is invertible; in particular, $s
\to 0$ iff $s' \to 0$. This means that the answers given by Theorem
3.4 for the $CC(Rj_{\nu*}\frak{F})$, glue together into one global
answer:
$$CC(Rj_*\frak{F}) = \lim_{s\to 0\atop s\in E}
(CC(\frak{F}) + \omega_s).\leqno (4.2)$$
Here the limit is taken in $T^*\Bbb P\,^n$ and $s$ runs over the set
$$
E = \{ (s_1,\dots,s_n) \in (\Bbb C^*)^n: \; s_i \neq s_j \; , \; i \neq j
\}.
$$
(The condition $s \in E$ is equivalent to $\phi_{\nu}^{-1}(s) \in
(\Bbb C^*)^{\{0,\dots,n\} \setminus \{\nu\}}$ for any $\nu$.)

Now, Theorem 1.3 for $G = (\Bbb C^*)^n$ would follow from (4.2) and the
next lemma.

\proclaim{(4.3) Lemma}
Let $\Lambda \subset T^*(\Bbb C^*)^n$ be an irreducible conic Lagrangian
variety. Then
$$\lim_{s\to 0\atop s\in E} (\Lambda +
\omega_s),[\Bbb P\,^n])_{T^*\Bbb P^n} = \roman{gdeg}(\Lambda).$$
\endproclaim

\noindent {\sl Proof:} Let $\Omega_s \subset T^*(\Bbb C^*)^n$ be the graph of
$\omega_s$. If $s \in E$, then $\Omega_s$ is closed in $T^*\Bbb P^n$ as
well. By translation, intersecting $\Lambda + \omega_s$ with
$[\Bbb P\,^n]$ is equivalent to intersecting $\Lambda$ with
$\Omega_{-s}$. By Proposition 2.2, there exists a Zariski open,
nonempty set $F \subset \Bbb C^n$ such that for $s^{(0)} \in F$ the
intersection $\Lambda \cap \Omega_{-s^{(0)}}$ consists of
$\text{gdeg}(\Lambda)$ smooth transverse points. Since $E$ is also Zariski
open, $F \cap E$ meets any polydisk
$$P_{\varepsilon} = \{(z_1,\dots,z_n) \in (\Bbb C^*)^n: \; 0 < |z_i| <
\varepsilon \}$$
in an open dense set. By the ``continuity of intersection'' it follows
that for $|\varepsilon | \ll 1$ and $s^{(0)} \in P_{\varepsilon} \cap F
\cap E$
$$
\biggl( \lim_{s\to 0\atop s\in E} \Lambda +
\omega_s),[\Bbb P\,^n]\biggl)_{T^*\Bbb P\,^n} = |\Omega_{-s^{(0)}} \cap \Lambda |
= \roman{gdeg}(\Lambda)
$$
and this completes the proof.

\heading{Proof of Theorem 1.3 in general}
\endheading

Let $A$ be an abelian variety. We recall some elementary properties of
line bundles on $A$, see [10] [16]. If $L$ is such a bundle, by $L_a$,
$a \in A$ we denote its fiber at $a$. By $L^0$ we denote the total
space of $L$ with the zero section deleted, so $L^0$ is a principal
$\Bbb C^*$-bundle on $A$. As for any base variety, the correspondence $L
\mapsto L^0$ is an equivalence between the category of line bundles
and isomorphisms and the category of principal $\Bbb C^*$-bundles.

Assume that $L$ has degree 0. Then the theorem of the square [16]
provides identifications
$$
L_0 \cong \Bbb C\; , \; L_a \otimes L_b \widetilde{\longrightarrow} L_{a+b}
$$
which make $L^0$ into a group, namely a semiabelian variety fitting
into an extension
$$0 \to \Bbb C^* \to L^0 \buildrel{p}\over\longrightarrow A \to 0.
\leqno (5.1)$$
Next, any bundle $L$ of degree 0 has a flat connection. All such
connections form an affine space $\text{Conn}(L)$ over the vector space
$H^0(A,\Omega^1) = \frak{a}^*$. Here $\frak{a}$ is the Lie algebra of
$A$. Let $\frak{L}$ be the Lie algebra of the group $L^0$, so that we
have the exact sequence
$$0 \to \frak{a}^* \to \frak{L}^* \buildrel{\pi}\over\longrightarrow
 \Bbb C\to
0.\leqno (5.2)
$$

Denote by $\widetilde{L} = p^*L$ the pullback of $L$ to $L^0$. For any
$\nabla \in \text{Conn}(L)$ let $\widetilde{\nabla}$ be its pullback to a
connection in $\widetilde{L}$. Denote also by $f$ the tautological section
of $\widetilde{L}$ over $L^0$ (given by the identity map).

\proclaim{(5.3) Proposition}
The 1-forms $\widetilde{\nabla} \log f$, $\nabla \in \roman{Conn}(L)$, are
invariant (with respect to the group structure on $L^0$). Their images
in $\frak{L}^*$ under the evaluation at 0 form the subspace
$\pi^{-1}(1)$, where $\pi$ is as in (5.2).
\endproclaim

\noindent{\sl Proof:} Let us first show the invariance. Note that
$\widetilde{\nabla} \log f$ is just the Lie algebra-valued 1-form on the
total space of the principal $\Bbb C^*$-bundle $L^0$, describing the connection
$\nabla$ in $L^0$ in the standard approach of differential
geometry. So its invariance follows from the fact that 
the line bundle $(L,\nabla)$
(or, what is equivalent, the $\Bbb C^*$-bundle $(L^0,\nabla)$) 
satisfies the theorem on the square as a bundle
with connection. More precisely, if $m,q_1,q_2: A \times A \to A$ are
the group structure and the two projections, then the isomorphism
$$
\mu: q_1^* L \otimes q_2^* L \to m^* L
$$
given by the theorem on the square is an isomorphism of bundles with
connection. In particular, the induced isomorphism $\mu_a: L_a \otimes
L \to l_a^* L$ is an isomorphism of bundles with connection on
$A$. But the translation $l_{(a,\lambda)}$, $\lambda \in L_a \setminus
\{ 0\}$ on $L^0$ is just given by $\mu_a (\lambda \otimes -)$. This
shows the invariance. As for the second assertion of Proposition 5.3,
it is again obvious from the interpretation of $\widetilde{\lambda} \log
f$ as the Lie algebra-valued 1-form describing the connection and the
identification $L_0^0 = \Bbb C^*$.

\proclaim{(5.4) Corollary}
The 1-forms $s \widetilde{\nabla} \log f$, $s \in \Bbb C^*$, $\nabla \in
\roman{Conn}(L)$, form a nonempty Zariski open set in the space of all
invariant 1-forms on $L^0$.
\endproclaim

Consider now several line bundles of degree 0, say $L_1,\dots,L_n$, on
$A$, and let
$$0 \to (\Bbb C^*)^n \to G = L_1^0 {\times}_A \dots
{\times}_A L_n^0\buildrel{p}\over\longrightarrow A \to 0
\leqno (5.5)$$
be the associated semiabelian variety.

We denote $\widetilde{L}_i = p^*L_i$ and let $f_i \in H^0(G,\widetilde L_i)$ be the
tautological section. As before, for $\nabla_i \in \text{Conn}(L_i)$ we
denote by $\widetilde{\nabla}_i$ its pullback to $\widetilde{L}_i$. Corollary
5.4 implies easily:

\proclaim{(5.6) Proposition}
Suppose $n \geq 0$. Then, the 1-forms
 $\sum _{i=1}^n s_i \widetilde{\nabla}_i \log f_i$ for
$s_1,\dots,s_n \in \Bbb C^*$, $\nabla_i \in \roman{Conn}(L_i)$, form a Zariski
open dense set in the space of all invariant 1-forms on $G$.
\endproclaim

We now turn to the proof of Theorem 1.3. 
Our approach is similar to [18], \S 2. First of all, it is known
that any semiabelian variety has the form (5.5) which we assume.
Next, we assume $n>0$ since for $n=0$ the group $G=A$ is compact and the
theorem follows from (1.1) and from Proposition 2.2. We set $L_0=\Cal O_A$ and compactify $G$ by embedding it into the relative projectivization
$$j: G\hookrightarrow \Bbb P := \Bbb P(L_0\oplus ... \oplus L_n)
\buildrel\rho\over\rightarrow A.$$
We can think of $\Bbb P$ as having homogeneous coordinates
$(t_0: ... ; t_n)$ with $t_i$ being not a function any more but rather
a section of $\rho^*L_i$. Note that $\rho^* L_i$ has a flat connection
induced from $\nabla_i$, whose restriction to $G$ is $\widetilde \nabla_i$. 
The variety $\Bbb P$ is the union of the
relative affine charts $A_\nu = \bigoplus_{i\neq \nu} L_i$. Inside $A_0$, 
the
complement of $G$ is given by
the condition that   one of the tautological sections (still denoted by $f_i$)
 of the pullback of $L_i$, vanishes. 

Take generic $\nabla_i\in\text{Conn}(L_i)$ so that for a
Zariski open, dense set of $s\in (\Bbb C^*)^n$ the 1-form
$\omega_s = \sum s_i\widetilde\nabla_i\log f_i$ satisfies
Proposition 2.2. 
Then, we mimic the arguments of Section 4 but in the relative
situation, using Theorem 3.5. We find that
$$CC(Rj_*\frak{F}) = \lim_{s\to 0\atop s\in E}
\biggl(CC(\frak{F}) + \sum_i s_i \widetilde{\nabla}_i \log f_i\biggr),
$$
the limit being taken in $T^*\Bbb P$. After this, the proof is
identical to the argument at the end of Section 4, and the theorem is
proven.

\Refs

[1]
A.A. Belinson, J.Bernstein, P. Deligne,  {\it Faiceaux
Pervers}, Asterisque {\bf 100} (1982), pp. 5-171.

\vskip .2cm

[2]
C.Chai, G. Faltings, Degeneration of Abelian
Varieties,  Springer-Verlag, 1990.

\vskip .2cm

{[3]}
J. Franecki, {The Gauss Map and Euler Characteristic on
Algebraic Groups}, Thesis,  Northwestern University, 1998.

\vskip .2cm

{[4]}
O. Gabber, F. Loeser,  {\it Faisceaux pervers $l$-adiques sur un
tore}, Duke Mathematical Journal {\bf 83} (1996), no. 3,
pp. 501-606. 

\vskip .2cm

{[5]}
I.M. Gelfand, M.M. Kapranov, A.V. Zelevinsky,  {\it Generalized Euler
integrals and $A$-hypergeometric functions}, Adv. in Math.
{\bf 84} (1990), no. 2, pp. 255-271.

\vskip .2cm

{[6]} I.M. Gelfand, M.M. Kapranov, A.V. Zelevinsky,
{\it  Hypergeometric functions,
toric varieties and newton polyhedra,} in: ``Special Functions"
( ICM-90 Sattelite Conference Proceedings,
M. Kashiwara, T. Miwa, Eds.) p. 104-121, 
Springer-Verlag, Berlin, 1991. 

\vskip .2cm

{[7]}
V. Ginsburg, {\it Characteristic varieties and vanishing cycles},
Invent. Math. {\bf 84} (1986), pp. 327 - 402.

\vskip .2cm

{[8]} V. Ginsburg, $\frak G$-{\it modules, Springer's representations
and bivariant Chern classes}, Adv. in Math. {\bf 61} (1986), 1-48.

\vskip .2cm

{[9]} V. Ginzburg, {\it Admissible modules on a symmetric space}, Ast\'erisque
{\bf 173-74} p.199-255, Soc. Math. France, 1989. 

\vskip .2cm

{[10]}
P. Griffiths, J. Harris, {Principles of Algebraic
Geometry}, John Wiley \& Sons, Inc., 1978.

\vskip .2cm

{[11]}
J. Horn, {\it \"{U}ber die Konvergenz hypergeometrischer Reihen zweier
und dreier Ver\"{a}nderlichen}, Math. Ann. {\bf 34} (1889),
pp. 544-600. 

\vskip .2cm

{[12]}
M.M. Kapranov,  {\it A characterization of $A$-discriminantal
hypersurfaces in terms of the logarithmic Gauss map}, Math. Ann. {\bf
290} (1991), pp. 277-285.

\vskip .2cm

{[13]}
M. Kashiwara,  {\it Index theorem for constructible sheaves},
Asterisque {\bf 130} (1985), pp. 193 - 209.

\vskip .2cm

{[14]}
M. Kashiwara, P.  Schapira, {Sheaves on Manifolds},
 Springer-Verlag, 1990.

\vskip .2cm

{[15]}
F. Loeser, C. Sabbah,  {\it Caract\'erisation des $D$-modules
hyperg\'eometriques irreductibles sur le tore}, Comptes Rendus
  Serie I. (Math\'ematiques), {\bf
312} (1991), no. 10, pp. 735 - 738, continued in {\bf 315} (1992),
no. 12, pp. 1263 - 1264. 

\vskip .2cm

{[16]}
D. Mumford, {Abelian Varieties}, Oxford University Press,
1970.
\vskip .2cm

{[17]}
W. Schmid, K. Vilonen, {\it Characters, characteristic cycles, and
nilpotent orbits}, in: {Geometry, Topology, and Physics (S.-T. Yau, Ed.)},
pp. 329-340, International Press, 1995.

\vskip .2cm

{[18]} P. Vojta, {\it Integer points on subvarieties of semiabelian varieties I,
} Invent. Math. {\bf 126} (1996), 133-181. 
\endRefs

\enddocument